\documentclass{article}
\usepackage{amssymb,amsmath,amsthm}
\usepackage{graphicx}

\def\bR{{\mathbb R}}

\def\cD{{\cal D}}

\def\Rq{{\bR^q}}

\def\Free{Free}

\newtheorem*{thmDH}{Theorem DH}

\newtheorem{example}{Example}

\newtheorem*{thm}{Theorem}
\newtheorem*{cor}{Corollary}

\def\cA{{\cal A}}
\def\Eq{g_{can}}

\title{A note on a conjecture of Gromov\\ about non-free isometric immersions}
\author{Roberto De Leo}
\begin{document}
\maketitle
\begin{abstract}
  We extend the results obtained in~\cite{DL07} towards the proof of a conjecture
  of M.~Gromov on isometric immersions via non-free maps. 
\end{abstract}
\section{Introduction}
The problem of isometric immersions was solved first by J.~Nash with his celebrated 
theorem~\cite{Nas56} and then, in a much wider setting, by M. Gromov~\cite{GR70,Gro86}.
In this context it was shown that, given a smooth $m$-dimensional manifold $M$, 
the operator
$\cD(f)=f^*\Eq$, which associates
to every smooth map $f:M\to\Rq$ its pull-back of the Euclidean metric $\Eq$ on $\bR^q$,
is an open map for $q\geq m(m+5)/2$ when restricted to the set of free maps
(throughout the paper we endow functional spaces with the Whitney 
topology).
We recall that free maps are maps
which are injective on the osculating space $J^2_0(\bR,M)$, i.e.
such that their first and second derivatives are linearly independent, and that
the set of all free maps $\mathop{Free}(M,\bR^q)$ is an open subset of 
$C^\infty(M,\bR^q)$ which is empty for $q < m(m+3)/2$ and dense for $q \geq m(m+5)/2$.

We focus here on the particular case $M=\bR^m$ and use $\cD_{m,q}$ for the operator 
$\cD$ acting on $C^\infty(\bR^m,\bR^q)$. In this case it is known that 
$\mathop{Free}(\bR^m,\bR^q)$ is non-empty for $q\geq m(m+3)/2$ (see Example~2 below) 
so that, in particular, it turns out that, for every $q\geq m(m+3)/2$, there is a 
non-empty open set $\cA$ on which the restriction of $\cD_{m,q}$ is an open map
(we say that such a $\cD$ is {\sl infinitesimally invertible} over the open set $\cA$);
very little instead is known about $\cD_{m,q}$ for smaller values of $q$, 
when no free map can arise for dimensional reasons. 
A conjecture by Gromov (see~\cite{Gro86}, p.162) claims that every $\cD_{m,q}$ 
is infinitesimally invertible over an open dense set for every $q \geq m(m+3)/2-\sqrt{m/2}$, 
namely that there exist non-empty dense open sets 
$\cA_{m,q}\subset C^\infty(\bR^m,\bR^q)$ such that $\cD_{m,q}|_{\cA_{m,q}}$ are open maps.

In a recent paper~\cite{DL07} steps were taken towards the full proof of the conjecture
by showing, through an explicit construction that made use of the Lie equations after
Gromov's idea in~\cite{Gro86} (p.~152), that $\cD_{2,4}$ is infinitesimally invertible 
over an open set $\cA^{DL}_{2,4}$.
In the next section we improve this result by finding a larger set $\cA_{2,4}\supset\cA^{DL}_{2,4}$ 
on which $\cD_{2,4}$ is open and defining open sets $\cA_{m,q}$, for every $q\geq m(m+3)/2-1$,
over which the operators $\cD_{m,q}$ are infinitesimally invertible. 
Our result is a direct consequence of
the following well-known theorem by Duistermaat and Hormander~\cite{DH72}:
%
\begin{thmDH}
  Let $M$ be an open manifold. Then the first order differential operator $X(f)=L_\xi f+\lambda f$, 
  where $\xi$ is a vector field on $M$ and $\lambda$ a smooth function, is surjective on 
  $C^\infty(M)$ iff $\xi$ admits a global transversal, i.e. an embedded hypersurface that cuts 
  in a single point each of its integral trajectories.
\end{thmDH}
\section{Proof and Examples}
Our aim is finding an open set of functions $\cA$ such that, if $f_0\in\cA$ and 
$g_0=\cD(f_0)$, the equation 
\begin{equation}
  \label{eq:D(f)=g}
  \cD(f)=g
\end{equation}
has solutions for every $g$ close enough to $g_0$.
We recall that, by a general theorem~\cite{Gro86}, the existence of solutions 
of~(\ref{eq:D(f)=g})
is granted by the existence of solutions of its linearized version, which in turn
is equivalent (e.g. see~\cite{GR70}) to the following algebraic system: 
\begin{equation}
  \label{eq:linsys}
  \begin{cases}
    \delta_{ij}\,\,\partial_\alpha f^i\delta f^j &=  h_\alpha\cr
    \delta_{ij}\,\,\partial_{\alpha\beta}f^i\,\delta f^j &= 
    (\partial_\alpha h_\beta + \partial_\beta h_\alpha - \delta g_{\alpha\beta})/2\cr
  \end{cases}
\end{equation}
where the $\delta f^i$ are the $q$ unknowns, 
the $\delta g_{\alpha\beta}$ are $m(m+1)/2$ given functions 
and the $ h_\alpha$ are $m$ arbitrary functions.
Therefore it is enough for our purposes to show that, for some open set of smooth functions, 
we can always choose the $ h_\alpha$ so that system~(\ref{eq:linsys}) has a solution.
\begin{thm}
  If $\displaystyle q\geq \frac{m(m+3)}{2}-1$ the 
  operators $\cD_{m,q}$ are infinitesimally invertible over non-empty open sets $\cA_{m,q}$.
\end{thm}
\begin{proof}
  When $q\geq m(m+3)/2$ the statement is trivially true because it is
  enough to choose $\cA_{m,q}=\Free(\bR^m,\bR^q)$.
  We will assume therefore in the remainder of the proof that $q=m(m+3)/2-1$, 
  i.e. that the number of equations is exactly one more than the number of unknowns $\delta f^i$. 

  Since the coefficients of the system~(\ref{eq:linsys}) are exactly the $m(m+3)/2$ 
  vector fields $\{\partial_\alpha f^i,\partial_{\alpha\beta}f^i\}$, then clearly there exist 
  non-identically zero functions $\lambda^\alpha$ and $\lambda^{\alpha\beta}=\lambda^{\beta\alpha}$, 
  such that identically
  $$
  \lambda^\alpha \partial_\alpha f^i + \lambda^{\alpha\beta} \partial_{\alpha\beta}f^i = 0.
  $$

  This reflects in the following compatibility condition for system~(\ref{eq:linsys}):
  $$
  2\lambda^\alpha  h_\alpha + \lambda^{\alpha\beta} 
  (\partial_\alpha h_\beta + \partial_\beta h_\alpha - \delta g_{\alpha\beta}) = 0  
  $$
  which is convenient to rewrite as 
  \begin{equation}
    \label{eq:ce}
    X^\alpha h_\alpha=\phi
  \end{equation}
  where $\phi=\lambda^{\alpha\beta}\delta g_{\alpha\beta}$, $X^\alpha$ is the first-order 
  non-homogeneous differential operator 
  $$
  X^\alpha=L_{\xi_\alpha}+2\lambda^\alpha,
  $$
  $L_{\xi_\alpha}$ is the Lie derivative with respect to the vector field 
  $$
  \xi_\alpha=\lambda^{\alpha\beta}\partial_\beta
  $$
  and the functions $\lambda^\alpha$ must be thought as the corresponding zero-order 
  multiplication operators.

  Now, let $\cA_{m,q}\subset C^\infty(\bR^m,\bR^q)$ be the open set of immersions $f$ such
  that their second-order jet $J^2_0(\bR,f):J^2_0(\bR,\bR^m)\to J^2_0(\bR,\Rq)$ has full rank
  and, in some coordinate system, there is a coordinate $\alpha_0$ of $\bR^m$ such that the 
  functions $\lambda^{\alpha_0 \beta}$ are never zero at the same time.
  Then, after setting $h_\beta=\lambda^{\alpha_0 \beta} h$, $\beta=1,\dots,m$, for some 
  unknown function $h$, the equation~(\ref{eq:ce}) becomes
  $$
  Y h=\psi
  $$
  where $Y=L_\zeta+\lambda'$ for some vector field $\zeta$ and function $\lambda'$. 
  A short computation shows that the component $\alpha_0$ of $\zeta$ is equal to
  $(\lambda^{\alpha_0 1})^2+\dots+(\lambda^{\alpha_0 m})^2$
  and therefore it is never zero by hypothesis. 
  In particular this means that every surface $x^{\alpha_0}=const$ is a global transversal 
  for $\zeta$ and therefore, by Theorem~DH, $Y$ is a surjective first-order partial
  differential operator. Hence for every function belonging to $\cA_{m,q}$ it is always 
  possible to choose the $h_\alpha$ in function of the $\delta g_{\alpha\beta}$ so 
  that the compatibility condition~(\ref{eq:ce}) is satisfied. 
\end{proof}
\begin{example}
  Consider any pair $(g,h)$ of free maps from $\bR$ to $\bR^2$. Then the function 
  $F_{gh}:\bR^2\to\bR^4$ defined by $F_{gh}(x,y)=(g(x),h(y))$ belongs to 
  $\cA_{2,4}\subset C^\infty(\bR^2,\bR^4)$. 
  Indeed in this case $\partial_{xy}F_{gh}=0$, so that we can choose 
  $$
  \lambda^x=\lambda^y=\lambda^{xx}=\lambda^{yy}=0\,,\,\,\lambda^{xy}=\lambda^{yx}=1
  $$ 
  and therefore the compatibility condition becomes simply
  $$
  \partial_x  h_y + \partial_y  h_x = \delta g_{xy}
  $$
  which is trivially solvable. 
  E.g. in concrete the function $F(x,y)=(x,e^x,y,e^y)$ belongs to $\cA_{2,4}$.
  Note that, while $\cA^{DL}_{2,4}\subset \cA_{2,4}$, $F$ does not belong to $\cA^{DL}_{2,4}$,
  i.e. $\cA_{2,4}$ is strictly bigger of the set introduced in~\cite{DL07}.
\end{example}
%
%
\begin{example}
  Let $F\in\Free(\bR^m,\bR^{m(m+3)/2})$ be the canonical free map given by
  $$
  F(x^1,\dots,x^m)=(x^1,\dots,x^m,(x^1)^2,x^1x^2,\cdots,(x^m)^2)
  $$
  and $\pi$ any projection $\pi:\bR^{m(m+3)/2}\to\bR^{m(m+3)/2-1}$ which ``forgets'' any 
  one of the last $m(m+1)/2$ components. Then their composition $F_\pi=\pi\circ F$ belongs
  to $\cA_{m,m(m+3)/2-1}$ because clearly the second-order jet of $F_\pi$ has full rank
  and one of its double derivatives, say $\partial_{x^1 x^2} F_\pi$, is identically zero, 
  so we can choose the corresponding factor $\lambda^{x^1 x^2}$ identically equal to 1.
  In the $(m,q)=(2,4)$ case for example we get the functions $F_1(x,y)=(x,y,xy,y^2)$,
  $F_2(x,y)=(x,y,x^2,y^2)$ and $F_3(x,y)=(x,y,x^2,xy)$.
\end{example}
Note that, exactly like in~\cite{DL07}, for $q=m(m+3)/2-1$ the set of $[m(m+3)/2]\times q$ matrices 
not satisfying the conditions introduced in the proof to define the open sets $\cA_{m,q}$ 
has just codimension 1 in the fibers of the bundle $J^2(\bR^m,\Rq)\to J^0(\bR^m,\Rq)$ while  
we would need at least codimension 3 in order to apply the transversality theorems.
In particular the sets $\cA_{m,q}$ are not dense.
\section*{Acknowledgments}
The author gladly thanks G.~D'ambra and A.~Loi for introducing the subject and for
many enlightning discussions and A.~Loi for a critical reading of the paper.
\bibliography{ii}
\end{document}